# Near-critical interior dynamics in a model of state–society interaction

Kerime Nur Kavadar[a], Ali Demirci[a,*] and Furkan Emre Isik[a]

[a]*Department of Mathematics, Faculty of Science and Letters, Istanbul Technical University, Istanbul, 34469, Turkey*

ABSTRACT

We study a normalized two–dimensional competitive Lotka–Volterra system describing the interaction between state power and society power. Restricting attention to the positively invariant domain $[0,1]^2$, the analysis focuses on interior–equilibrium dynamics where coexistence persists as the unique long–run outcome. We show that when the interaction parameters approach the coexistence threshold $a_{12}a_{21} \to 1^-$, convergence toward equilibrium becomes slow and trajectories exhibit prolonged transient dynamics. In this near–critical regime, trajectories organize into a narrow corridor around the balance manifold despite the absence of bistability. The corridor structure can be characterized quantitatively through equilibrium gaps and interaction thresholds. Numerical simulations illustrate how asymmetric adjustment influences the geometry and persistence of these transient regimes.

## 1. Introduction

Competitive Lotka–Volterra systems originate from the seminal work of Volterra [1] and constitute a classical class of nonlinear dynamical models for studying interaction, dominance, and coexistence under mutual inhibition. Their qualitative theory has been extensively developed in the dynamical systems literature, particularly for low-dimensional competitive flows [2, 3], while normalized formulations are commonly employed to analyze bounded interaction dynamics [4].

Beyond biological applications, competitive dynamics have also been used in political economy to describe the interaction between state capacity and societal power. A central reference is Acemoglu and Robinson [5], who distinguish weak, despotic, and inclusive states and emphasize that inclusive regimes—characterized by mutual constraint and coexistence—are both desirable and intrinsically fragile. Related studies highlight how institutional capacity, political incentives, and collective constraints jointly shape long-run outcomes [6, 7, 8].

In much of this literature, fragility is associated with regime transitions or bistability. However, continuous-time competitive systems suggest a different mechanism: fragility may arise within a single coexistence regime through weakening attraction and prolonged transient dynamics. In such situations, the equilibrium remains unique and locally stable, yet trajectories may spend long periods near intermediate states before convergence occurs.

Motivated by this perspective, we study a normalized two-dimensional competitive Lotka–Volterra system defined on a positively invariant unit square. The model admits a unique interior equilibrium whenever the coexistence condition $a_{12}a_{21} < 1$ holds. Our analysis focuses on the near-critical regime in which the interaction parameters approach the threshold $a_{12}a_{21} \to 1^-$. In this limit the determinant of the Jacobian becomes small, producing a pronounced separation of time scales: trajectories rapidly contract toward a low-dimensional attracting direction and subsequently evolve slowly toward equilibrium.

This geometric mechanism generates extended transient phases in which the two state variables remain close to each other over long time intervals. We refer to this transient region of the phase space as a *narrow corridor*. The phenomenon provides a continuous-time dynamical interpretation of fragile balance in state–society interaction: coexistence persists as the unique asymptotic outcome, yet convergence may become arbitrarily slow near the coexistence threshold. Numerical illustrations show how interaction asymmetries and adjustment rates influence the persistence and structure of these corridor-like transients.

---

*Corresponding author
✉ kavadar20@itu.edu.tr (K.N. Kavadar); demircial@itu.edu.tr (A. Demirci); isikf17@itu.edu.tr (F.E. Isik)
ORCID(s): 0000-0002-9144-5691 (A. Demirci)

## 2. Model and preliminaries

We consider a two-dimensional competitive dynamical system describing the interaction between two aggregate quantities, denoted by $L(\tau)$ (state power) and $S(\tau)$ (society power), where $\tau$ represents dimensionless time. Throughout this study, both variables are assumed to be *normalized*, so that

$$0 \leq L(\tau) \leq 1, \qquad 0 \leq S(\tau) \leq 1, \tag{1}$$

for all admissible times $\tau \geq 0$.

The normalization reflects a modeling perspective in which the state variables represent relative levels or intensities rather than absolute magnitudes. Such formulations are particularly suitable when the focus lies on bounded interaction dynamics and comparative dominance under implicit capacity constraints [4]. Accordingly, the admissible phase space is restricted to the compact domain

$$\Omega := [0,1] \times [0,1]. \tag{2}$$

This choice ensures that the dynamics remain confined to a finite and interpretable state space.

Within this normalized framework, we study the competitive system

$$\begin{aligned}
\frac{dL}{d\tau} &= L(1 - L - a_{12}S), \\
\frac{dS}{d\tau} &= \rho S(1 - S - a_{21}L),
\end{aligned} \tag{3}$$

where $a_{12}, a_{21} > 0$ denote cross-interaction coefficients. The parameter $\rho > 0$ represents the relative adjustment rate of $S$ with respect to $L$, introducing an asymmetric time-scale structure into the system. Models of this form arise naturally in competitive Lotka–Volterra dynamics with self-limitation and asymmetric adaptation mechanisms [9].

A fundamental property of system (3) is that the domain $\Omega$ is positively invariant. On the boundary $L = 0$ or $S = 0$, the corresponding right-hand side vanishes, preventing trajectories from leaving $\Omega$ through the coordinate axes. Moreover, on the boundaries $L = 1$ and $S = 1$, the vector field points inward. Consequently, solutions initiated in $\Omega$ remain in $\Omega$ for all future times, and $\Omega$ constitutes the natural phase space for the analysis.

System (3) admits three boundary equilibrium points,

$$(0,0), \quad (1,0), \quad (0,1), \tag{4}$$

as well as a potential interior equilibrium

$$E^* = (L^*, S^*), \tag{5}$$

provided that

$$a_{12}a_{21} < 1. \tag{6}$$

In this case, the interior equilibrium is given by

$$L^* = \frac{1 - a_{12}}{1 - a_{12}a_{21}}, \qquad S^* = \frac{1 - a_{21}}{1 - a_{12}a_{21}}. \tag{7}$$

The interior equilibrium $E^*$ lies in $\Omega^\circ$ provided that $0 < a_{12} < 1$ and $0 < a_{21} < 1$, which ensure that both $L^*$ and $S^*$ are strictly positive.

Although competitive Lotka–Volterra systems may, in a global sense, exhibit bistable dynamics associated with competing boundary attractors [10], such regimes are excluded in the present framework. The normalized domain $\Omega$ is introduced precisely to focus on interior dynamics, where both variables remain strictly positive and bounded away from degenerate states. Boundary equilibria therefore correspond to limiting configurations and are not regarded as admissible long-term outcomes.

Under the above parameter restrictions, the interior equilibrium $E^*$ governs the asymptotic behavior of trajectories initiated in $\Omega^\circ = (0,1)x(0,1)$. The Jacobian matrix evaluated at $E^*$ has negative trace and positive determinant, implying local asymptotic stability [9]. In particular,

$$\det J(E^*) = \rho L^* S^* (1 - a_{12}a_{21}). \tag{8}$$

To quantify proximity to the coexistence threshold, we introduce the small parameter

$$\eta := 1 - a_{12}a_{21}. \tag{9}$$

The regime $\eta \ll 1$ corresponds to interaction parameters close to the coexistence boundary. In this limit the determinant of the Jacobian becomes small and one eigenvalue of the linearization approaches zero in magnitude, indicating the emergence of slow dynamics near equilibrium.

The asymmetric time-scale parameter $\rho$ does not affect the location of equilibrium points, but it influences transient behavior by controlling the relative speed at which $S$ adjusts compared to $L$. This effect becomes particularly significant in the near-critical regime and will play an important role in the transient dynamics analyzed in the next section.

## 3. Near-critical dynamics and corridor formation

We now analyze the transient behavior of system (3) in the near-critical regime

$$\eta = 1 - a_{12}a_{21} \ll 1.$$

Although the interior equilibrium $E^*$ remains locally asymptotically stable for all $\eta > 0$, the approach to equilibrium changes qualitatively as $\eta$ becomes small.

The Jacobian matrix at the interior equilibrium is

$$J(E^*) = \begin{pmatrix} -L^* & -a_{12}L^* \\ -\rho a_{21}S^* & -\rho S^* \end{pmatrix}.$$

Its trace and determinant are given by

$$\operatorname{tr} J(E^*) = -(L^* + \rho S^*) < 0, \qquad \det J(E^*) = \rho L^* S^* \eta.$$

Hence the eigenvalues $\lambda_{1,2}$ satisfy

$$\lambda_1 + \lambda_2 = -(L^* + \rho S^*), \qquad \lambda_1 \lambda_2 = \rho L^* S^* \eta.$$

When $\eta \to 0^+$, one eigenvalue remains $O(1)$ while the other scales as $O(\eta)$. This produces a clear pronounced separation of time scales in the local dynamics near $E^*$.

The spectral structure described above implies that trajectories approaching the equilibrium exhibit two distinct phases: a rapid contraction toward a one-dimensional attracting direction, followed by slow evolution along this direction toward equilibrium.

**Proposition 3.1.** *Let $E^*$ denote the interior equilibrium of system (3). For $\eta = 1 - a_{12}a_{21} \ll 1$, trajectories starting in $\Omega^\circ$ approach a neighborhood of $E^*$ along a fast direction and subsequently converge to equilibrium on a slow time scale of order $O(\eta^{-1})$.*

This behavior is characteristic of near-critical competitive systems, where weakening hyperbolicity generates long transient phases despite the persistence of a unique stable equilibrium.

The separation of time scales leads to a geometric organization of trajectories in the phase space. After the initial fast contraction, solutions remain confined for long time intervals within a narrow region in which the two state variables remain close to each other before the final approach to equilibrium.

We refer to this transient region as a *narrow corridor*. Importantly, the corridor does not correspond to an invariant set or an alternative attractor. Rather, it emerges from the near-critical geometry of the flow and the slow eigen-direction associated with the small eigenvalue.

Consequently, even though coexistence remains the unique asymptotic outcome, trajectories may spend long periods near balanced states when $\eta$ is small. This mechanism provides a dynamical explanation for prolonged phases of relative balance observed in the numerical simulations presented in Section 5.

## 4. Regime interpretation and political economy implications

In this section we interpret the qualitative regimes of system (3) through the lens of political economy. Although the analysis is conducted within a normalized dynamical framework, the equilibrium configurations admit natural interpretations in terms of institutional balance and dominance.

System (3) admits three boundary equilibria, $(0,0)$, $(1,0)$, and $(0,1)$, in addition to a potential interior equilibrium $E^* = (L^*, S^*)$. From a mathematical perspective these correspond to degenerate or limiting steady states. In a political economy interpretation, however, they represent extreme institutional configurations.

The equilibrium $(0,0)$ corresponds to a collapse of both components, which may be interpreted as a breakdown of institutional capacity or legitimacy on both sides of the interaction. The boundary equilibrium $(1,0)$ represents unilateral dominance of the state, while $(0,1)$ corresponds to the dominance of society. Such corner outcomes are often associated with exclusionary or extractive institutional arrangements in which balance and mutual constraint are absent.

In the present framework these boundary configurations are not regarded as relevant long-term outcomes. The normalized domain $\Omega$ is introduced precisely to focus on interior dynamics in which both actors remain active. Consequently, the analysis concentrates on the interior equilibrium $E^*$, which represents a regime of coexistence between state and society.

Importantly, coexistence does not imply symmetry. The equilibrium levels $L^*$ and $S^*$ depend on the interaction coefficients $a_{12}$ and $a_{21}$, reflecting asymmetric influence between the two actors. Different parameter values therefore correspond to regimes that are state-leaning or society-leaning, even though both components remain strictly positive.

The most interesting dynamics arise near the coexistence threshold. As shown in Section 3, when the interaction parameters approach the critical condition $a_{12}a_{21} \to 1^-$ the system becomes near-critical: one eigenvalue of the linearized dynamics becomes small and trajectories exhibit a pronounced separation of time scales. Solutions rapidly approach a low-dimensional attracting direction and then evolve slowly toward equilibrium.

This geometric mechanism produces extended transient phases during which the two components remain close to each other before the final approach to equilibrium. We interpret this region of the phase space as a *narrow corridor* of balanced interaction. Within this corridor, coexistence is preserved but dynamically fragile: small perturbations may lead to prolonged deviations before the system returns to equilibrium.

From a political economy perspective, this regime corresponds to a situation in which balance between state and society is maintained through continuous adjustment rather than through strong restoring forces. The resulting dynamics capture the idea that inclusive institutions may remain stable in the long run while still exhibiting persistent and fragile short- and medium-term fluctuations.

## 5. Numerical illustrations

In this section we present numerical simulations illustrating the transient dynamics of the normalized competitive system (3). All simulations are performed within the positively invariant domain $\Omega = [0,1]^2$ using nonnegative initial conditions. Parameter values are chosen so that a unique interior equilibrium exists in each case, in agreement with the analytical results. The numerical experiments are intended to visualize how the interaction parameters influence transient behavior, particularly in the near-critical regimes described in Section 3.

All trajectories are computed using the MATLAB `ode45` solver with relative and absolute tolerances set to `RelTol=`$10^{-6}$ and `AbsTol=`$10^{-8}$. To visualize the corridor behavior in the time domain, we introduce a simple indicator based on the difference between the two state variables. Let $(L^*, S^*)$ denote the interior equilibrium. We define the equilibrium-based tolerance

$$\varepsilon := |L^* - S^*|, \tag{10}$$

which measures the intrinsic asymmetry between the equilibrium components.

Using this quantity, we introduce the normalized deviation

$$\gamma(\tau) = \frac{|L(\tau) - S(\tau)|}{\varepsilon}. \tag{11}$$

Time intervals for which $\gamma(\tau) \leq 1$ correspond to phases in which the two variables remain close relative to the equilibrium gap. This indicator provides a convenient numerical way to visualize the corridor region predicted by the near-critical dynamics discussed in Section 3.

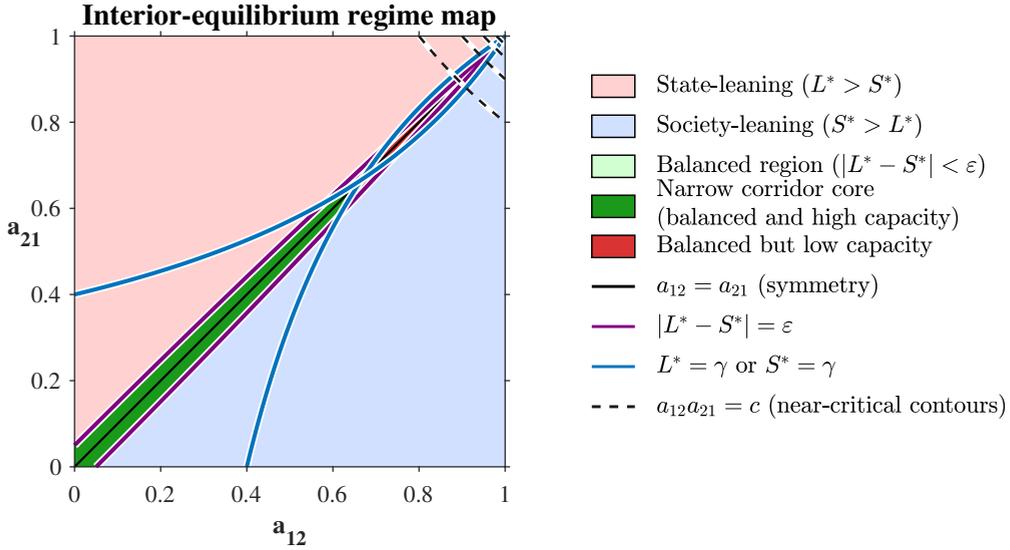

**Figure 1:** Interior–equilibrium map on the interaction–parameter plane $(a_{12}, a_{21}) \in [0,1]^2$, constructed from the explicit interior equilibrium $E^* = (L^*, S^*)$ of system (3). Colors indicate equilibrium-based classifications: state-leaning regimes ($L^* > S^*$), society-leaning regimes ($S^* > L^*$), and a balanced region defined by $|L^* - S^*| < \varepsilon$, where $\varepsilon$ denotes the equilibrium gap used in the corridor indicator. The narrow-corridor core corresponds to balanced and high-capacity configurations satisfying $|L^* - S^*| < \varepsilon$ together with $\min\{L^*, S^*\} > \gamma$, while the balanced low-capacity region satisfies $|L^* - S^*| < \varepsilon$ and $\min\{L^*, S^*\} \leq \gamma$. Overlaid curves indicate the symmetry line $a_{12} = a_{21}$, the equilibrium-gap contour $|L^* - S^*| = \varepsilon$, and the capacity contours $L^* = \gamma$ and $S^* = \gamma$, which delineate transitions between the corresponding regimes.

We first summarize the equilibrium-level structure of the system on the interaction-parameter plane $(a_{12}, a_{21}) \in [0,1]^2$. Figure 1 presents the interior equilibrium map obtained from the explicit expressions

$$L^* = \frac{1 - a_{12}}{1 - a_{12}a_{21}}, \qquad S^* = \frac{1 - a_{21}}{1 - a_{12}a_{21}},$$

whenever $a_{12}a_{21} < 1$.

The classification is based on two equilibrium indicators: the gap $|L^* - S^*|$ and the minimum equilibrium level $\min\{L^*, S^*\}$. A balanced region is identified by $|L^* - S^*| < \varepsilon$, while a high-capacity balance condition is imposed by $\min\{L^*, S^*\} > \gamma$, where $\varepsilon > 0$ and $0 < \gamma < 1$ are chosen thresholds used in the numerical construction. The symmetry line $a_{12} = a_{21}$ and the curves $|L^* - S^*| = \varepsilon$ together with $L^* = \gamma$ or $S^* = \gamma$ are overlaid to indicate transitions between the corresponding regimes.

We next illustrate representative time-domain behavior in Figure 2. First consider parameters close to the coexistence threshold,

$$a_{12} = 0.48, \qquad a_{21} = 0.55, \qquad \rho = 1.$$

For this configuration the interior equilibrium lies strictly inside $\Omega$, while the interaction parameters remain close enough to the coexistence boundary to produce extended transient phases. The corresponding time series show long intervals in which $\gamma(\tau) \leq 1$, indicating that $L(\tau)$ and $S(\tau)$ remain close to each other before the final convergence to equilibrium. This behavior is consistent with the slow convergence mechanism associated with the near-critical spectral structure described in Section 3.

We also examine parameter configurations that favor one component over the other. For

$$a_{12} = 0.50, \qquad a_{21} = 0.70, \qquad \rho = 1,$$

the transient dynamics display temporary dominance of $L(\tau)$ before converging to the interior equilibrium. Conversely, when

$$a_{12} = 0.75, \qquad a_{21} = 0.45, \qquad \rho = 1,$$

the trajectories exhibit transient dominance of $S(\tau)$.

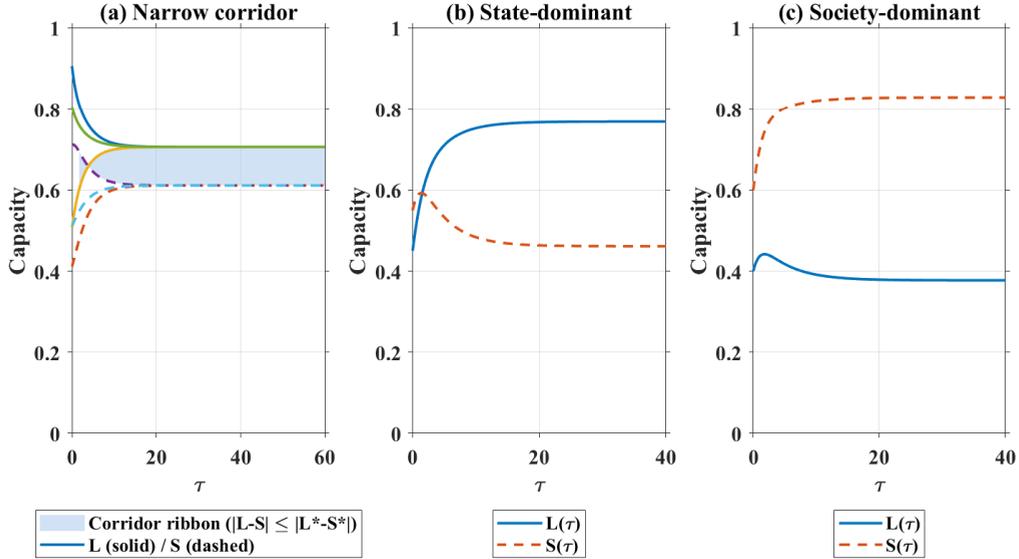

**Figure 2:** Transient time-series dynamics of system (3) under different interaction structures with $\rho = 1$. Panel (a) corresponds to parameters $a_{12} = 0.48$ and $a_{21} = 0.55$ with three different initial conditions, which lie close to the coexistence threshold. The trajectories exhibit extended intervals in which the corridor indicator satisfies $\gamma(\tau) \leq 1$, indicating that $L(\tau)$ and $S(\tau)$ remain close for long periods before convergence. Panels (b) and (c) illustrate asymmetric interaction regimes. For $a_{12} = 0.50$, $a_{21} = 0.70$ the transient dynamics display temporary dominance of $L(\tau)$, whereas for $a_{12} = 0.75$, $a_{21} = 0.45$ the trajectories exhibit transient dominance of $S(\tau)$. In all cases the long-term behavior is governed by the same interior equilibrium, while differences arise primarily in the structure and duration of the transient phase, consistent with the near-critical dynamics discussed in Section 3.

These examples demonstrate that variations in the interaction coefficients primarily affect the transient structure of the dynamics. In all cases the long-term behavior remains governed by the same interior equilibrium, while the duration and character of the transient phase depend on the proximity to the coexistence threshold and the degree of interaction asymmetry.

## Conclusion

We studied a normalized competitive dynamical system with asymmetric adjustment, focusing on interior-equilibrium dynamics within the positively invariant domain $\Omega = [0, 1]^2$. As the interaction parameters approach the coexistence threshold $a_{12}a_{21} \to 1^-$, the equilibrium remains locally stable but convergence becomes slow. This near-critical regime generates prolonged transient dynamics organized around a narrow region of approximate balance. The resulting corridor-like structure arises without bistability: coexistence remains the unique asymptotic outcome while trajectories remain close to balance for extended periods. Numerical simulations illustrate how asymmetric adjustment influences the persistence and geometry of these transient corridors.

From a political economy perspective, the model provides a continuous-time interpretation of the narrow corridor idea in which balance between competing forces can persist despite remaining dynamically fragile. Within this framework the corridor can be characterized quantitatively through equilibrium gaps and interaction thresholds, offering a mathematical criterion for identifying when balance becomes fragile yet persistent. The formulation also suggests natural extensions in which additional interacting components represent a third center of influence in the balance of power, or where interaction parameters evolve over time. In this sense, simple competitive dynamical systems may offer a tractable mathematical framework for quantifying and extending the narrow corridor concept in models of institutional balance.

## References

[1] V. Volterra, Fluctuations in the abundance of a species considered mathematically, Nature 118 (1926) 558–560.
[2] M. Hirsch, S. Smale, R. Devaney, Differential Equations, Dynamical Systems, and an Introduction to Chaos (2nd ed.), Academic Press, 2004.
[3] I. M. Bomze, Lotka–volterra equation and replicator dynamics: a two-dimensional classification, Biological Cybernetics 48 (1983) 201–211.
[4] J. D. Murray, Mathematical Biology I: An Introduction, Springer, 2002.
[5] D. Acemoglu, J. A. Robinson, The emergence of weak, despotic and inclusive states, National Bureau of Economic Research (2017) No:w23657.
[6] T. Besley, T. Persson, The origins of state capacity: Property rights, taxation, and politics, American Economic Review 99 (2009) 1218–1244.
[7] T. Besley, T. Persson, State capacity, conflict, and development, Econometrica 78 (2010) 1–34.
[8] D. Acemoglu, Politics and economics in weak and strong states, Journal of Monetary Economics 52 (2005) 1199–1226.
[9] L. Perko, Differential Equations and Dynamical Systems (3rd ed.), Springer, 2006.
[10] J. Hofbauer, K. Sigmund, Evolutionary Games and Population Dynamics, Cambridge University Press, 1998.